\documentclass{IEEEtran4PSCC}
\usepackage{amsmath,amssymb, amsfonts,amsthm}
\usepackage{algorithm,algorithmic,bbm}

\usepackage{hyperref}
\newtheorem{theorem}{Theorem}

\newtheorem{lemma}[theorem]{Lemma}

\newtheorem{remark}[theorem]{Remark}

\newtheorem{DE}{Definition}[section]
\newtheorem{PRO}[DE]{Procedure}
\newtheorem{TEMP}[DE]{Template}

\newcommand{\mA}{{\mathcal{A}}}
\newcommand{\mB}{{\mathcal{B}}}
\newcommand{\mK}{{\mathcal{K}}}
\newcommand{\mF}{{\mathcal{F}}}
\newcommand{\mE}{{\mathcal{E}}}
\newcommand{\mG}{{\mathcal{G}}}

\newcommand{\mR}{{\mathcal{R}}}
\newcommand{\mS}{{\mathcal{S}}}

\newcommand{\mW}{{\mathcal{W}}}
\newcommand{\bmath}[1]{{\mathbf{#1}}}
\newcommand{\bE}{{\mathbf{E}}}
\newcommand{\bP}{{\mathbf{P}}}
\newcommand{\bT}{{\mathbf{T}}}
\newcommand{\bVar}{{\mathbf{Var}}}
\newcommand{\bV}{{\mathbf{V}}}
\newcommand{\bStd}{{\mathbf{Std}}}
\usepackage{ifpdf}

\usepackage{cite}

  \usepackage[pdftex]{graphicx}
\usepackage{array}
\usepackage{fixltx2e}

\usepackage{stfloats}
 \usepackage{url}

\begin{document}
%
\title{Variance-Aware Optimal Power Flow}

\author{
\IEEEauthorblockN{Daniel Bienstock and Apurv Shukla}
\IEEEauthorblockA{
Columbia University, NY,
United States\\
\{dano, apurv.shukla\}@columbia.edu}
}

\maketitle

\begin{abstract}
  The incorporation of stochastic loads and generation into
  the operation of power grids gives
  rise to an exposure to stochastic risk.  This risk has been addressed in
  prior work through a variety of mechanisms, such as scenario generation or
  chance constraints, that can be incorporated into OPF computations.  Nevertheless, numerical experiments reveal that the
  resulting operational decisions can produce power flows with very
  high variance.  In this paper we introduce a variety of convex 
  variants of OPF that explicitly address the interplay of (power flow) variance
  with cost minimization, and present numerical experiments that highlight
  our contributions.
\end{abstract}


\noindent 

\section{Introduction}
Increasing levels of renewable penetration have resulted in associated risk,
principally through potential equipment overloads.  A growing
body of research has tackled this issue through proposed modifications to
power engineering practices, in particular by suggesting alternative, risk-aware formulations to the Optimal Power Flow (OPF) problem that may additionally
incorporate appropriate
balancing mechanisms.  See e.g. \cite{line0, maria0, ccopf, miles, maria1, line2, morrison, phan, line1, zhangli, sundar2016unit}, \cite{miles}   and citations therein.  

The application of these methodologies yields operational plans within the
scope of an OPF computation that guarantee protection from stochastic
variations in generation (or loads). For example, in the chance-constrained
case, one obtains a set of dispatch decisions and possibly participation
factors (which describe balancing control) that guarantee that the probability that any given line becomes
overloaded is at most a given value $\epsilon > 0$, under appropriate assumptions on renewable stochasticity. 

However, experiments
with such schemes produce cases where power flows are highly
variable, or more precisely some transmission lines will have power flows
with high \textit{variance}.  We posit that
this behavior is undesirable from a real-time operational perspective.  See \cite{hiskensvariability}, which highlights the impact of power flow variability
on voltage profile and transformer operation. High variability
may also hinder system understanding
and control, 
as well as pricing mechanisms (see \cite{davari}, \cite{jabr1}, \cite{yamin} for
work on pricing under uncertainty).  \cite{schellenberg} considers minimization of variance of slack bus injection. 

In this paper we propose a set different alternatives to risk-aware OPF computations,
which we term \textit{variance-aware OPF}.  In these variants we explicitly
account for power flow variance, either as a post-processing step for the
initial risk-aware OPF computation, or by explicitly incorporating
a measure of power flow variance into a risk-aware OPF formulation so as to
obtain a convex optimization problem.

\section{Notation and basic formulations}\label{notaform}
In this paper we will focus on the linearized, or DC model for power flows. We
will use the following selected nomenclature, with additional terms described
later:\\

\noindent $\mB$ = set of buses, $n = |\mB|$; $B$ = bus susceptance matrix.\\
\noindent $\hat{B}=B$, with the last row and column removed \\
\noindent $\breve{B}=\left[{\begin{array}{cc} \hat{B}^{-1} &0 \\ 0 & 0 \end{array}}\right] $; $\breve{B}_i$ = $i^{th}$ row of $\breve{B}$. \\
\noindent $\theta_i$ = phase angle at bus $i$, $\bar \theta_i \ = \ \bE(\theta_i)$.\\
\noindent $\mE$ = set of lines; $m = |\mE|$. \\
\noindent $f^{\max}_{ij}$, $b_{ij}$ =  power flow limit, and susceptance, for line $ij \in \mE$.\\
\noindent $f_{ij}$ =  power flow on line $ij$, $\bar f_{ij} = \bE(f_{ij}) = b_{ij} (\bar \theta_i - \bar \theta_j)$  \\
\noindent $\pi_{ij} \, \doteq \, \breve{B}^T_i - \breve{B}^T_j$ for each line $ij \in \mE$\\
\noindent $\mG$ = set of generator buses; we assume at most one generator\\
 \hspace*{.2in} per bus\\
 \noindent $p^{min}_i, \, p^{max}_i$ = minimum and maximum output of $i \in \mG$.\\
\noindent $\mS$ = set of stochastic injection buses, 
\noindent $\mu_k \, + \, \omega_k$ = stochastic injection at bus $b \in \mB$,
\begin{itemize}
\item $\mu_k$ = constant, $\omega_k$ = zero-mean random variable.
   \item $\mu_k = \omega_k$ = 0 if $k \notin S$,
   \item $\Omega$ = covariance of $\omega$, an $n \times n$ matrix,
   \item $\mW$ = set of vectors $\omega$ under consideration.
\end{itemize}
\noindent $\mR$ = set of buses participating in balancing\\
\noindent $\mA$ = matrix of participation factors; $\alpha_{ij} =$ (i,j)-entry of $\mA$,
\begin{itemize}
  \item $\mA$ is $n \times n$
  \item $\alpha_{ij} = 0$ if $i \notin \mR$ or $j \notin \mS$
\item $\bV(\mA)$ = vector with entries $b_{ij}^2 \pi_{ij}^T(I - \mA) \Omega (I - \mA^T) \pi_{ij}$ (variances) for $ij \in \mE$
\end{itemize}
\noindent $\mK$ = convex set of allowable participation factor matrices \\
\noindent $p$ = vector of generation amounts, extended to all $b \in \mB$\\
\hspace*{.1in}(forcing $p_b = 0$ when $b \notin \mB$).\\
\noindent $c(p)$ = (convex) cost of generation vector $p$.\\
\noindent $d$= (fixed) vector of loads.\\

\noindent Using this notation, we can write a DC-OPF formulation: 
\begin{subequations}
\label{stdOPF}
\begin{eqnarray}
  && \min_{p} c(p) \\
\mbox{s.t.}   && B\theta \ = \ p \, - \, d \label{balance}\\ 
&& \forall ij \in \mE: \ b_{ij}|\theta_i - \theta_j| \, \leq f_{ij}^{\max} \label{linelimit} \\
&& \forall i \in \mG: \ p_{i}^{min} \leq p_{i} \leq p_{i}^{max} \label{genlimit}
\end{eqnarray}
\end{subequations}
As is well-known, constraints \eqref{balance}-\eqref{linelimit} can be simplified by using a pseudo-inverse for the matrix $B$ (see, e.g. \cite{line0}, \cite{ccopf}).  Namely, for each
line $ij$ the power flow on $ij$ equals $b_{ij} \pi^T_{ij} (p - d)$ thereby reducing \eqref{balance}-\eqref{linelimit} to the system
\begin{eqnarray}
  && \forall ij \in \mE: b_{ij} | \pi^T_{ij} (p - d) | \, \leq f_{ij}^{max} \label{linelimit2}
\end{eqnarray}
\section{Security-constrained formulations}
Previous work \cite{line0}, \cite{ccopf}, \cite{miles}, \cite{line1}, \cite{line3}, \cite{morrison}, \cite{maria0}, \cite{maria1}   has focused on modifications
to power flow computations in order to account for uncertain injections, and
in particular to model the use of balancing. We will first outline some of this work.

Second, we will also describe a sparse formulation for security-constrained DC-OPF problems, as a straightforward extension of the approach used in \cite{line0}, \cite{ccopf}, \cite{line1}. This will be the starting point for the discussions in Section \ref{formulations}, where we will discuss how to make  the formulation
variance-aware.

The method used in \cite{line0}, \cite{ccopf}, \cite{line1},
modifies the DC-OPF computation so that its primary output
consists of a vector $\bar p$ of \textit{controllable generation} amounts
(which can only be positive at controllable generators) and (as in \cite{line1})
an $n \times n$ matrix $\mA$ used to model balancing. Given $\omega \in \mW$ the balancing vector is given by $\mA \omega$, so that the stochastic output of a controllable generator at bus $i$ is
\begin{align}
& p_i \, = \, p_i(\omega) \, = \, \bar p_i \ - \ [\mA \omega]_i  \ = \ \bar p_i \ - \ \sum_{j \in \mB} \alpha_{ij} \omega_j \label{controllable}
\end{align}
where we assume that when $i$ is not a controllable
generator bus then $\alpha_{ij}= 0$ for all $j$. Thus, the net injection vector equals $ \bar p - d + \mu + \omega - \mA \omega$. 
In order to actually attain balancing we must have
\begin{eqnarray}
\forall w \in \mW:  && \sum_{i \in \mB}(\bar p - d + \mu + \omega - \mA \omega)_i \ = \ 0, \label{balanceguarantee}
\end{eqnarray}
a stochastic requirement. Let us assume that we additionally assume $\sum_{i \in \mB}(\bar p - d + \mu )_i = 0$, i.e. the system is balanced when $\omega = 0$.  Then \eqref{balanceguarantee} is equivalent to
\begin{eqnarray}
\forall \omega \in \mW:  && \omega_i \, - \, \sum_{j \in \mB} \alpha_{ij} \omega_j \ = \ 0, \ \ \forall \, i \in \mB.  \label{shorter}
\end{eqnarray}
As pointed out in \cite{line1} in order
to attain this condition it is sufficient to require 
\begin{eqnarray}
  && 1 \ = \ \sum_{j \in \mB} \alpha_{ij} \quad \forall \, i \in \mB. \label{supershort}
\end{eqnarray}
Note that \eqref{shorter} describes a hyperplane in $\omega$-space.  If $\mW$ is
full dimensional then \eqref{supershort} is actually necessary for \eqref{balanceguarantee} \cite{batt}.  Additionally we may constrain $\mA$ in a number of ways
(e.g. the $\alpha_{ij}$ constrained by bounds).  In particular, we could 
require that for any bus $i \in \mR$ \footnote{Henceforth, a participating bus.}, $\alpha_{ij} = \alpha_{ik}$ for every $j \neq k$
(termed a ``global'' policy in \cite{line1}).

\noindent {\bf Notation.} In what follows we will write $\bmath{\mA \in \mK}$ to denote a set of generic convex constraints that $\mA$ must satisfy,
including in particular \eqref{supershort}.

Using these notations, if $f_{ij}$ denotes (stochastic) flow on line $ij$ then
\begin{subequations} \label{Eexpression}
\begin{eqnarray}
\hspace*{-.2in}\forall ij \in \mE:  && f_{ij} \, = \, b_{ij} \pi^T_{ij} (\bar p- d + \mu + \omega - \mA \omega), \\
\hspace*{-.2in}&& \bE(f_{ij}) \, = \, b_{ij} \pi^T_{ij} (\bar p- d + \mu). 
\end{eqnarray}
\end{subequations}
Likewise by construction
\begin{align}
  & \bV(\mA)_{ij} \, \doteq \ \bVar(f_{ij}) =b_{ij}^2 \bVar(\pi_{ij}^T (I - \mA)\omega) = \nonumber\\
  &   \quad b_{ij}^2 \pi_{ij}^T(I - \mA) \Omega (I - \mA^T) \pi_{ij} \label{Varexpression}
\end{align}

The next modeling ingredient concerns security constraints. A
variety of  variants of the
basic problem \eqref{stdOPF} can be obtained depending on
how we model stochasticity and security and on whether the matrix $\mA$ is an optimization variable or is fixed (i.e. given as an input). A typical approach concerns
\textit{chance constraints}.  To fix ideas, consider a line $ij$.  Then we wish
to impose that 
\begin{align} \label{chanceline}
& \bP( | f_{ij}|\, > \, f_{ij}^{\max}) \, < \, \epsilon
\end{align}
where $ 0 < \epsilon < 1$ is a given threshold.  In order to represent \eqref{chanceline} in a convex manner, prior work has assumed normally distributed of $\omega$.  Under such an assumption \eqref{chanceline} is equivalent to
\begin{align} \label{secondchancegeneric}
& | \bE (f_{ij})| \, + \,   \ \Phi^{-1}(1 - \epsilon) \bStd (f_{ij}) \ \le \ f^{\max}_{ij}
\end{align}
where $\Phi^{-1}(1 - \epsilon)$ is the $\epsilon$-quantile for a normal distribution.  Using \eqref{Eexpression} and \eqref{Varexpression}, we obtain that \eqref{secondchancegeneric} is SOCP representable \cite{line0}, \cite{ccopf} (also see below).

\subsection{Modifications used in this paper}
Our first modification will be to replace \eqref{secondchancegeneric} with
\begin{align} \label{safety}
& | \bE (f_{ij})| \, + \,   \ \nu_{ij} \, \bStd (f_{ij}) \ \le \ f^{\max}_{ij}
\end{align}
where $\nu_{ij}$ is a \textit{safety parameter}.  As seen above, in
the Gaussian case \eqref{safety} is
equivalent to \eqref{chanceline} if we choose $\nu_{ij} = \Phi^{-1}(1 - \epsilon)$.
However, the Gaussian case is not the only one
where such an equivalence holds; other examples of distributions of interest  include (multivariate) truncated Gaussians, uniform distributions on
ellipsoidal supports, and others.  This requires
a distribution-dependent choice $\nu_{ij} = \nu_{ij}(\epsilon)$.
Further \eqref{safety} can be used
to tightly approximate distributionally robust chance constraints.  See \cite{calafioreelghaoui} (and e.g. Theorem 2.1 therein).   Similar remarks apply to security
constraints involving generators. We would argue that even when the stochastics of $\omega$ is complex so that
we cannot provide a rigorous choice for the safety parameters, we may still be able to
compute 
$\Omega$, or perhaps a data-driven estimate for it.  The safety-parameter
approach could still hold appeal from an intuitive, if imprecise, perspective.  We will term \eqref{safety} a \textit{safety} constraint.

In order to present our modification to chance-constrained DC-OPF, we first produce
a new expression for the variance of a line flow $f_{ij}$. Let us
write
\begin{eqnarray}
D & \doteq & \breve B \mA, \quad \mbox{and} \\
\forall \, k \in \mS:   \  \gamma_{ij,k} & \doteq & \breve B_{ik} - \breve B_{jk} - D_{ik} + D_{jk}.
\end{eqnarray}
\begin{lemma} \label{varlemma}
  For any line $ij$, the variance of flow on $ij$ under scheme \eqref{controllable} is given by
  \begin{eqnarray} \label{newVarexpression}
    \bVar(f_{ij}) & = & b_{ij}^2 \, \gamma_{ij} \, \Omega \, \gamma_{ij}^T
  \end{eqnarray}
\end{lemma}    
\noindent {\bf Remark:}  we stress that this expression holds without any assumption on the underlying probability distributions. Denote the generation cost at a bus $i$ as
$c_i(p) \doteq c_{i0}p^2 + c_{i1} p + c_{i0}$, and let the $i^{th}$ row of  $\mA$ be denoted by $\mA_i$.  Then (routine proof):
\begin{lemma} \label{costlemma} For a bus $i$, $\bVar(p_i) = \mA_i^T \Omega \mA_i$, and 
$\bE( c_i(p_i)) = c_{i0}(\bar p_i^2 + \mA_i^T \Omega \mA_i) + c_{i1} \bar p_i + c_{i2}.$
\end{lemma}

Our safety-constrained formulation, given next, generalizes the chance-constrained
formulation in \cite{ccopf}.
As before we write $n = |\mB|$ and also $m = |\mE|$.
We are given nonnegative a safety parameter $\nu_{ij}$ for each line $ij$ and
likewise safety parameters  $\nu_i$ for generators ($= 0$ at non-generator buses).
The formulation uses variables
$\bar p, \bar \theta$ ($n$-vectors), $\bar f$ ($m$-vector), $\mA, D$ ($|\mS| \times |\mS|$ and $n \times |\mS|$ matrices, respectively),
 and $\gamma$ and $s$ (an $m \times |\mS|$ matrix and $m$-vector, respectively).  As above, we use $\mA \in \mK$ to denote a given set of convex constraints on $\mA$.
\begin{subequations}
\label{scOPF1}
\begin{align}
\hspace*{-.1in} & \min \quad \sum_{i \in \mG} \bE( c_i(p_i)) \label{scOPFobj} \\
\hspace*{-.1in} & \text{s.t.} \ \mA \in \mK \label{genAK}\\
\hspace*{-.1in}   & B\bar{\theta}=\bar{p}+\mu-d \label{flowxpression} \\
\hspace*{-.1in}   & \bar f_{ij} \ = b_{ij}(\theta_i - \bar \theta_j) \\
\hspace*{-.1in}   & b_{ij} | \bar f_{ij}| + \nu_{ij} \, s_{ij} \, \le \, f^{\max}_{ij} \quad \forall ij \in \mE, \label{tail}\\
\hspace*{-.1in}   & \breve B \mA \ = D \label{Ddef}\\
\hspace*{-.1in}   & \gamma_{ij, k} \ = \ \breve B_{i,k} - \Breve B_{j,k} - D_{i,k} + D_{j,k}, \quad \forall ij \in \mE, \, k \in \mB \label{gammadef}\\
\hspace*{-.1in}   & s_{ij} \, \ge \, b_{ij} \sqrt{ \gamma_{ij} \Omega \gamma_{ij}^T} \quad \forall ij \in \mE\label{theconic}\\
\hspace*{-.1in} & \forall i \in \mG: \nonumber \\
\hspace*{-.1in} & \quad p_{i}^{min} + \nu_i \sqrt{ \mA_i^T \Omega \mA_i} \ \leq \bar p_{i} \leq p_{i}^{max} - \nu_i \sqrt{ \mA_i^T \Omega \mA_i}. \label{genlimit2}
\end{align}
\end{subequations}
Equations \eqref{Eexpression} and Lemma \eqref{varlemma}, together with constraints
\eqref{flowxpression}-\eqref{theconic} imply that \eqref{tail} correctly states the desired 
line safety constraint.  Similarly with the generator safety constraints.  Problem \eqref{scOPF1} is a convex quadratically constrained
optimization problem that can, in principle, be solved using standard optimization software.  We call this the \textit{sparse} formulation.  Next we analyze the structure of this formulation, in light of the fact that
previous work \cite{line1}, \cite{ccopf}, \cite{line3} has highlighted numerical difficulties in solving chance-constrained formulations for relatively large systems. Thus it is important to point out that our formulation in fact can be greatly
reduced in size.

In particular, as denoted in Section \ref{notaform} and in the remarks following
Lemma \ref{varlemma} if a bus $j \notin S$ then $\alpha_{ij} = 0$ for every
bus $i$,  and if bus $i \notin \mR$
then $\alpha_{ij} = 0$ for every $j$.  In summary we have:
\begin{lemma} \label{sparsesizelemma}
  In the sparse formulation, the number of $\mA$, $D$ and $\gamma$ variables is, respectively, $ |\mR| |\mS|$, $n |\mS|$ and $ m \mS|$.  In addition we also have the $\bar \theta$, $\bar f$ and $p$ variables, which number $2 n + m$. The number of nonzeros in  all constraints \eqref{Ddef} and \eqref{theconic} is $ n |\mR| |\mS|$ and $ m |\mS|$, respectively.
\end{lemma}

We note here are several variants of the above safety-constrained problem that could be
meaningful.  First, there has been recent work on how to avoid the
normality assumption \cite{tillmann}.  Another alternative would be to rely
on scenario modeling (see e.g. \cite{maria0}), or to use
a distributionally robust model with underlying normality \cite{miles}. 

\section{A numerical example}\label{exp1}
In this section we consider some simple examples of security-constrained problems \eqref{scOPF1} so as to examine the tradeoff between cost and various measures of line flow variance.  Consider the following example:
\begin{figure}[h]
\centering
\includegraphics[scale=.37]{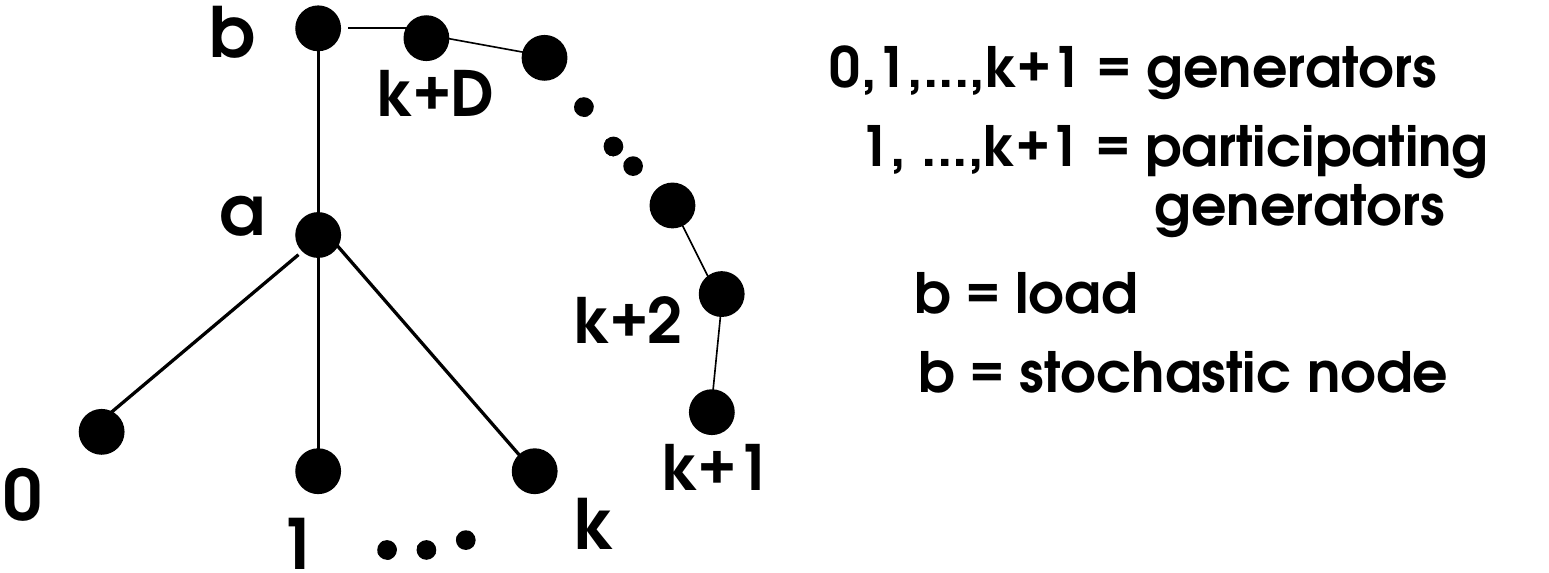}
\caption{High-variance example.}
\label{fig:1}
\end{figure}

Here,
\begin{itemize}
\item Quantities $k$ and $D$ are large 
\item Bus $b$ has a load of $L$ units.
\item The stochastic output at bus $b$  is indicated by $\omega$, with mean
  $\mu < L$ and variance $\sigma^2$.
\item Bus $i$ ($0 \le i \le k+1$) is a generator with linear cost
  function given by $c_{i1}p_i$.\\ Assume $c_{01} < c_{11} =c_{21} \ldots c_{k1} < c_{k+1,1}$.
\item The generator at bus $0$ is large (capacity larger than $L$) and non-participating.
\item  Generators at buses $1,2, \ldots, k+1$ are all participating, with
  zero lower limit and safety parameters of value $3$.
\item Buses $a, k+1, \ldots, k+D$ have no load, no generation and no stochastic
  injection.
\item We assume (for the first analysis) that all line limits are large.
\end{itemize}
Let us write the participating factor for generator $i \ge 1$ as $\alpha_i$.  
Consider the following candidate solution:
\begin{subequations}\label{candidate}
\begin{align}
& \text{$\bar p_0 = L - \mu - 3 \sigma$}\\
& \text{$\alpha_i = 1/k$ and $\bar p_i = 3 \sigma/k$, for $1 \le i \le k$.  }\\
& \text{$\alpha_{k+1} = \bar p_{k+1} = 0$}.
\end{align}
\end{subequations}
\begin{lemma}\label{foo}
  The solution stipulated by \eqref{candidate} is the unique optimal solution for the case of problem \eqref{scOPF1} in Figure \ref{fig:1}.
\end{lemma}
\noindent \textit{Proof.} First we argue that the solution we gave is feasible. The generation at bus $i$ ($1 \le i \le k$) equals
$3\sigma/k - \omega/k$. Hence the security constraint at bus $i$ is satisfied.
Also, total generation equals $L - 3 \sigma + 3 \sigma - \omega = L$.  Thus
indeed the solution is feasible.
Now consider any other feasible solution  given by values $\hat p_i$ and $\alpha_i$ (average
generation and participating factor at bus $i$).  Define $P = \sum_{i = 1}^{k+1} \hat p_i$.  Then $\hat p_{0} = L - \mu - P$, and so the cost of the solution is
\begin{align}
& c_{01} (L - \mu - P) + c_{11} (P - p_{k+1}) + c_{k+1,1} \hat p_{k+1} \ge \\
&  c_{01} (L - \mu - P) + c_{11} P.
\end{align}
Here the inequality is strict unless $p_{k+1} = 0$.  Further, $\hat p_{i} \ge 3 \alpha_i \sigma$ and $\sum_i {\alpha_i} = 1$. So, $P \ge 3 \sigma$.  The proof is
now complete. \qed

Let us now consider the stochastic flow on line $ab$, which equals $L - \mu - \omega$, and therefore has variance $\sigma^2$.  In other words line $ab$ is exposed to $100 \%$ of stochastic injection variance.  (Additional comments, below).  The sum of line flow variances equals $\sigma^2 (1 + 1/k) \approx \sigma^2$ for  $k$ large.

Suppose we were to aim for a decrease of variance on this line by
$50 \%$.  This goal will be achieved by setting $\sum_{i = 1}^k \alpha_i = \sqrt{.5}$ and thus $\alpha_{k+1} \approx 0.293$.  In that case the sum of variances will be
larger than $ .5 \sigma^2 + (D+1) \alpha_{k+1}^2 \sigma^2 \approx (.5 + (.293)^2 (D+1))\sigma^2$.  With $D = 10$ this quantity equals approximately $1.44 \sigma^2$ and the sum of line flow variances has substantially increased!\\

In the previous example we assumed that line limits were large.  However, it is
simple to adapt the example to that were line limits are small enough that
the safety constraints become significant.  At the same time, it is of interest to consider an alternative variance metric that
takes into account line limits. A compelling metric is the sum,
over all lines, of the ratio of flow variance  to square of line limit.
To address this metric we modify the data in the example above by assuming:
\begin{itemize}
\item  $\mu = L/4$ and $\sigma = \mu/2 = L/8$.
\item The limit for lines $0a$ and $ab$ equals $9L/8$.
\item The limit for all lines $ia$ ($1 \le i \le k$)  and for the lines on
  the path from $k+1$ to $b$ are all equal to $2 \sigma$.
\item All line safety parameters have value $3$.  
\end{itemize}
We can verify that the solution given by \eqref{candidate} remains feasible,
in which case a close reading of Lemma \ref{foo} shows that it remains the
sole optimal solution.   To verify feasibility, note that the expected flow
on line $ab$ equals $3L/4$ and its standard deviation equals $L/8$; since
the limit for that line is $9L/8$ it follows that the safety constraint is
satisfied (exactly), and likewise with the remaining lines.  Thus feasibility of
\eqref{candidate} holds, and it is straightforward to show that Lemma \ref{foo}
still holds.

Further, let us consider the variance metric (sum of ratios of line flow variances to square of line limits) in solution \eqref{candidate}.  Line $ab$ contributes $1/81$ to the sum whereas each line $ia$ contributes equals $1/(4k^2)$.  Hence the total metric equals $1/81 + 1/(4 k) \approx .0123$ for $k$ large.  If, as
above we desire to decrease the variance on line $ab$ by $50 \%$ then 
we again set $\alpha_{k+1} = 1 - \sqrt{.5} \approx .293$.  We also set $\bar p_{k+1} = 3 \alpha_{k+1}$.  Since the line limit on the path from $k+1$ to $b$ is $2 \sigma > 6 \alpha_{k+1}$ it follows that this solution is indeed feasible for \eqref{scOPF1}.  The
variance metric for the new solution is (more than) $0.5/81 + (.293)^2 (D+1)/4 \approx 0.242$. In other words, in our attempt to reduce variance in one line
we have worsened the variance metric by approximately a factor of $20$. If
we reduce $D$ from $10$ to $3$ the metric for the new solution is approximately
$0.098$ and thus more than $8$ times larger than that for the original solution.\\

To conclude this section we point out that:
\begin{enumerate}
\item In either case, a cursory analysis might seem to reveal that any solution to the
  security-constrained problem will concentrate variance on line $ab$.  That is not so: as we have shown, we can shift variance to the path from bus $k+1$ to $b$.  This approach, however,
  both increases cost and system variance.  Also note that on the surface the network design above makes sense: cheaper participating generators are closer to
  the load.
\item   The system we have discussed has been expressly designed so as to evince an extreme behavior.  And, of course, different parameter choices will result in less extreme behavior. However we caution the reader that the example can be modified by introducing more participating generators and additional loads, with the appropriate cost and security structure, to obtain a system setup that is not at first glance extreme, while the actual behavior does become even more extreme than the one
we provided.
Moreover the qualitative behavior
we have discussed \textit{can be seen in realistic examples}.
\item The behavior we see in the example, in summary, is one where network topology, load structure (i.e. location of loads, stochastic injection nodes and responding buses) and cost structure conspire so as to concentrate variance on a specific line or more generally a small set of lines.  
\end{enumerate}
As the example shows, we can expect an inherent set of tradeoffs between operational cost, system-wide variance metrics, and variance on selected (``important'') lines.  This tradeoff is one that should be visible to and controllable by system operators.  Next section addresses this point.
\section{Variance-aware problems}\label{formulations}
In this section we consider modifications to the process of solving formulation \eqref{scOPF1}, so as to capture
tradeoffs with variance metrics.  As a \textit{starting point} in this direction, we formulate an optimization problem of the general type
\begin{subequations}
\label{varOPFcommon}
\begin{align}
  \hspace*{-.1in} & \min \quad \sum_{i \in \mG} \bE( c_i(p_i)) \ + \ \Delta(\bar f, s^2) \label{varOPFcommonobj} \\
\hspace*{-.1in} & \text{s.t.} \ \text{\ref{genAK} - \ref{genlimit2}}
\end{align}
\end{subequations}
where $s^2$ is the vector with entries $s^2_{ij}$ (variances of line flows), and
\begin{align}
  & \Delta(\bar f, s^2)  \ = \ \sum_{ij \in \mF}  \Delta_{ij}(\bar f_{ij}, s^2_{ij}), \label{Deltadef}
\end{align}  
Here, $\mF \subseteq \mE$ and each $\Delta_{ij}$ is a nonnegative function chosen to highlight a specific
penalty as a function of expected flow and variance.  In other words, we impose all
constraints of the safety-constrained problem \eqref{scOPF1} but we add to the
objective a term that enforces a tradeoff with variance. The function $\Delta$ will be our
formal ``variance metric.''
\subsection{Some variance metrics}\label{varmetric}
Next we describe some concrete models of variance metrics \eqref{Deltadef}.
\begin{itemize}
\item [(I)] $\mF = \mE$, and for all $ij \in \mE$, $\Delta_{ij}$ is convex and nondecreasing in $s^2_{ij}$. As a special case, let $\Delta_{ij}(\bar f_{ij}, s_{ij}) = \psi_{ij} s_{ij}^2$ where $\psi_{ij} \ge 0$.  When $\psi_{ij} = 1$ or   $\psi_{ij} = (1/f_{ij}^{max})^2$ $\Delta(\bar f, s^2)$ is the metric we considered in the examples given in Section \ref{exp1} . 
\item [(II)] More complex models are those where $\mF = \mF(\bar f, s^2)$.  For example let $N > 0$ be given, and consider a function $\Delta$ of the form
\begin{align}
  \hspace*{-.1in} & \Delta(\bar f, s) \ = \ \sum_{ij \in \mF} \psi_{ij} s^2_{ij} 
\end{align}
where again the $\psi_{ij} \ge 0$ are scalars, and $\mF \subseteq \mE$ is a set of lines, such as\\
\hspace*{.2in} (II.1) The set of $N$ lines with largest flow magnitude.\\
\hspace*{.2in} (II.2) The  set of $N$ lines with largest flow variance.\\
In either model, we do not know in advance the set $\mF$ to be summed over, yet the
objective is relevant.  
\item [(III)] A very different class of convex models is that where for each $ij$ we
  define $\Delta_{ij}(\bar f_{ij}, s_{ij})$ as follows:
  \begin{align}
\Delta_{ij} = \   & -\rho_{ij} \log( s^2_{ij} - b^2_{ij} \gamma_{ij} \Omega \gamma_{ij}^T) ,\quad \text{if} \ s^2_{ij} > b^2_{ij} \gamma_{ij} \Omega \gamma_{ij}^T \nonumber\\
=  \ & +\infty, \quad \text{otherwise}. \nonumber
\end{align}    
where $\rho_{ij} > 0$ is a parameter.
This is the classical \textit{logarithmic barrier} formulation \cite{nocedalwright}.  Because of the definition of the $\Delta_{ij}$, the conic constraints \eqref{theconic} can be removed. It is known that if the $\rho_{ij}$ are all equal to a common value $\rho$, then the solution to 
\eqref{varOPFcommon} (with \eqref{theconic} removed) converges to an optimal solution to the security-constrained problem \eqref{scOPF1} as $\rho \rightarrow 0^+$. For a fixed choice
of the $\rho_{ij}$ we obtain a ranking of the importance of (security of) individual lines, plus a tradeoff against generation cost.  The formulation \eqref{varOPFcommon} may thus be seen as a viable
alternative to safety-constrained formulations, without the 
computational burden of the conic constraints \eqref{theconic}.
\end{itemize}
In model (I) we obtain a standard second-order cone program.  Model (III) is also convex, though nonstandard (i.e. not an SOCP).  Model (II.2) can be formulated as a convex program (proof omitted).

\subsection{A correction template}
Our goal in this work is to address the tradeoff between generation cost, security constraints and variance metrics, and one way to do so is to solve the
optimization problem \eqref{varOPFcommon}, which is in principle a straightforward task in model (I).  However, previous experience with chance-constrained DC-OPF \cite{ccopf}, \cite{line2} indicates that a direct solution approach relying
on state-of-the-art solvers is likely to fail due to numerical difficulties. Our numerical tests with
formulation \eqref{varOPFcommon} in model (I) verifies this fact.

As an alternative to formal optimization we instead focus on a procedure that seeks to \textit{correct} or \textit{adjust} the solution to the non-variance aware safety-constrained problem \eqref{scOPF1}.  A template for the overall scheme is as follows.

\begin{center}
  \fbox{
    \begin{minipage}{0.9\linewidth}
      \hspace*{.9in} \begin{TEMP}{GENERIC CORRECTION TEMPLATE}\label{gencasc} \end{TEMP}
              {\bf Input}: an instance of the safety-constrained problem \eqref{scOPF1} and a variance metric.\\ \\
              {\bf Step I.} Solve \eqref{scOPF1}, with solution $(\bar p^*, \mA^*)$.\\ \\
              {\bf Step II.} Perform a small number of adjustment iterations which shift $(\bar p^*, \mA^*)$ to a new feasible solution to \eqref{scOPF1} that attains an improved value of the variance metric, while at the same time increasing generation cost in a moderate manner.
    \end{minipage}
  }
\end{center}
We will describe an implementation of Step {\bf II} where we perform first- or second-order steps that amount to solving convex optimization problems which do not include large numbers of conic constraints.

Prior to describing the specific implementation we briefly discuss the motivation for using this template.  As has been observed in prior work, in the solution to typical chance-constrained problems similar to \eqref{scOPF1},
only very a small number of the conic constraints \eqref{theconic} are binding or
nearly binding, that is to say, only a small number of safety constraints are
active.  If some of those lines are among those with largest flow (literally, a handful) we expect that their variance will also be fairly small.  However, there generally is a
relatively large number of lines with quite large variance of flow.  It follows
that there is ``low hanging fruit'', i.e. opportunities for shifting variance so as to improve a given variance metric.  In doing so, we may slightly increase
the  variance of lines with nearly-binding safety-constraints.  This
increase may require a correspondingly slight decrease of flow in such lines, so as to satisfy
the safety constraints.  Flow will thus be shifted onto other lines,
possibly resulting in a small increase in cost. 
These expectations are borne out in experiments detailed below, and they form the underpinning for the above template.

\subsection{A specific implementation for the template}\label{tempimp}
Here we describe an implementation of the above template. {\bf For simplicity of exposition, we assume that}
\begin{itemize}
\item[{\bf (1)}] For each line $ij$, $\Delta_{ij}(\bar f_{ij}, s^2_{ij}) = \Delta_{ij}(s^2_{ij})$, i.e. it is a function of line variance only.  Further, it is assumed that $\Delta_{ij}(s^2_{ij})$ is convex and
  nondecreasing.
\item [{\bf (2)}] The set $\mF$ used to define $\Delta(\bar f, s^2)$   in \eqref{Deltadef} does not depend on the variances, i.e. $\mF = \mF(\bar f)$.\\
  
  Assumptions (1) and (2) match models (I) and (II.1) above. However we stress that
  there are variants of the procedures below that handle cases (II.2) and (III). 
\end{itemize}

\begin{DE} Let $\bar f$ be a power flow vector and let $\mA \in K$ be a
matrix of participation factors.
We will say that the pair $(\bar f, \mA)$ is
{\bf compatible} (or that $\mA$ is compatible with $\bar f$) for problem \eqref{scOPF1}, if there exist $\bar p, \bar \theta, D, \gamma$ and $s$ such that $(\bar p, \bar f, \bar \theta, \mA, D, \gamma, s)$ is feasible for \eqref{scOPF1}.
\end{DE}
Informally, a pair $(\bar f, \mA)$ is compatible if they give rise to a feasible solution
to the safety-constrained problem \eqref{scOPF1}.

Our implementation relies on two optimization problems that are repeatedly solved in Procedure \ref{varshift} given below.  We describe these two problems next. Let $0 < \tau < 1$ be fixed. The first optimization
problem is denoted by $\bmath{Reroute(\hat \mA, \tau)}$, and uses as inputs a compatible pair $(\hat f, \hat \mA)$ and vector $\hat s^2 = \bV(\mA)$ of line flow variances arising from  $\hat \mA$. 
\begin{subequations}
\label{scOPFfixedvar}
\begin{align}
\hspace*{-.1in} & \min_{\bar p, \bar f, \bar \theta} \quad \sum_{i \in \mG} c_{i0}(\bar p_i^2 + \hat \mA_i^T \Omega \hat \mA_i) + c_{i1} \bar p_i + c_{i2} \label{scOPFobjfixed} \\
\hspace*{-.1in} & \text{s.t.} \quad B\bar{\theta}=\bar{p}+\mu-d \label{flowxpressionfixed} \\
\hspace*{-.1in}   & b_{ij} | \bar f_{ij}| + \nu_{ij} \, \hat s_{ij} \, \le \, (1 - \tau) f^{\max}_{ij} \quad \forall ij \in \mE, \label{tailfixed}\\
\hspace*{-.1in}   & \bar f_{ij} \ = \ b_{ij} (\bar \theta_i - \bar \theta_j) \quad \forall ij \in \mE, \\
\hspace*{-.1in} & \forall i \in \mG: \nonumber \\
\hspace*{-.1in} & p_{i}^{min} + \nu_i \sqrt{ \hat \mA_i^T \Omega \hat \mA_i} \ \leq \bar p_{i} \leq p_{i}^{max} - \nu_i \sqrt{ \hat \mA_i^T \Omega \hat \mA_i}. \label{genlimit2fixed}
\end{align}
\end{subequations}
\noindent {\bf Comments:} This optimization problem minimizes expected generation cost using  
the fixed participation factors $\hat \mA$; it imposes stricter line safety constraints (with line limits reduced by the factor $1 - \tau$).  Assuming $\bmath{Reroute(\hat \mA, \tau)}$ is feasible, let an optimal solution  be $(\bar p^*, \bar f^*, \bar \theta^*)$.  Then, by construction of $\bmath{Reroute(\hat \mA, \tau)}$,
$(\bar f^*, \hat \mA)$ is compatible for \eqref{scOPF1}, with some slack. A large choice for $\tau$ may of course render $\bmath{Reroute(\hat \mA, \tau)}$ infeasible. However we will be choosing small values for $\tau$, which as a corollary implies that
the expected generation cost accrued by $\bar p^*$ will be slightly larger than that entailed by the flow $\hat f$.\\
The second problem takes as input a compatible pair $(\bar f', \mA')$.
Also, for $ij \in \mE$ let $s'_{ij} = \sqrt{\bV(\mA')_{ij}}$, the standard deviation of flow on $ij$
using participation factors $\mA'$. Finally, let $\bT(\bar f', \mA', \tau)$ be the set of lines for which the safety constraint is
nearly tight under participation factors $\mA'$, that is to say:
$$ \bT(\bar f', \mA', \tau) \ = \ \{ \, ij \in \mE \, : \, |\bar f'_{ij}| + \nu_{ij} s'_{ij} \ge (1 - \tau)f^{max}_{ij} \}.$$
The problem, 
denoted by $\bmath{VShift(\bar f', \mA', \tau)}$, is as follows.  
\begin{subequations}
\label{scOPFshiftvar}
\begin{align}
\hspace*{-.1in} & \min_{s, \mA} \quad \sum_{ij \in \mF(f')}\Delta_{ij}(s_{ij}^2) \\
\hspace*{-.1in} & \text{s.t.} \quad \mA \in \mK\\
\hspace*{-.1in}   & s^2_{ij} \ \ge \ b_{ij}^2 \pi_{ij}^T(I - \mA) \Omega (I - \mA^T) \pi_{ij} \quad \forall \, ij \in \mE \label{sdef}\\
\hspace*{-.1in}   & | \bar f'_{ij}| \, + \,   \ \nu_{ij} \, s_{ij} \ \le \ f^{\max}_{ij} \quad \forall ij \in \bT(f', \mA', \tau). \label{careful}
\end{align}
\end{subequations}
\noindent {\bf Comments.}  Since the $\Delta_{ij}$ are increasing, at
optimality constraint \eqref{sdef} will be binding. Thus the optimization problem is selecting a participation matrix that minimizes the variance metric.
Constraint \eqref{careful} stipulates that on lines $ij \in \bT(f', \mA', \tau)$
the matrix $\mA$ is compatible with the input flow $f'$.  Problem \eqref{scOPFshiftvar} is convex and when the $\Delta_{ij}$ are quadratic it is an SOCP; its relative difficulty depends on the number of constraints \eqref{careful} which, as we have discussed, is frequently quite small.\\

We now use the two optimization problems to develop an algorithm.  

\begin{center}
  \fbox{
    \begin{minipage}{0.9\linewidth}
      \hspace*{.9in} \begin{PRO}{Variance-shifting}\label{varshift} \end{PRO}
              {\bf Input}: Feasible solution $\bmath{(\bar p^0, f^0, \mA^0)}$ to safety-constrained problem \eqref{scOPF1}, variance metric $\bmath{\Delta}$, parameters $0 < \tau < 1$, $K > 0$.   Let  $\bmath{s^2_0 = V(\mA_0)}$.\\ \\
      {\bf For} $k = 1, 2, \ldots, K \ $ {\bf perform iteration k:} \\ \\
      \hspace*{.05in}{\bf 1.} Solve $\bmath{Reroute(\mA_{k-1}, \tau)}$.\\
\hspace*{.25in}{\bf If} infeasible, {\bf STOP.}\\
      \hspace*{.25in}{\bf Else,} let $\bmath{(\bar p^k, \bar f^k, \bar \theta^k)}$ be the optimal solution.\\ \\
      \hspace*{.05in}{\bf 2.} Solve $\bmath{VShift(\bar f^k, \mA_{k-1}, \tau)}$, with solution $\bmath{(\hat s_k, \hat \mA_k)}$.\\ \\
        \hspace*{.05in}{\bf 3.} Choose $0 < \lambda \le 1$ largest, so that \\
        \hspace*{.2in} $\bmath{(\bar f^k, (1 - \lambda) \mA_{k-1} + \lambda \hat \mA_k)}$ is compatible. \\ \\
        \hspace*{.05in}{\bf 4.} Set $\bmath{ \mA_k \ \leftarrow \ (1 - \lambda) \mA_{k-1} + \lambda \hat \mA_k}$,\ $\bmath{s^2_k = V(\mA_k)}$.\\ \\
        \hspace*{.05in}{\bf 5.} If $\bmath{\Delta( \bar f^k, s^2_k) \ge \Delta(\bar f^{k-1}, s^2_{k-1})}$. {\bf STOP}.\\
        \hspace*{.2in} Reset $\tau \leftarrow \tau/2$.
    \end{minipage}
  }
\end{center}

A formal analysis of this algorithm, which is motivated by Nesterov's smoothing techniques for non-smooth problems \cite{nesterov}\footnote{With $\Delta$ playing the role of a ``potential'' function.}, is provided in the Appendix.  Intuitively, in Step 1 the algorithm reroutes flow so as so as to create slack capacity in lines, while keeping a constant participation matrix (so that variances remain constant).  In Step 2 we
compute a new participation matrix \textit{which is not required} to be compatible
with the flow vector computed in Step 1, but should improve on the variance
metric because of the typically small number of constraints \eqref{careful}.

Finally, in Steps 3-4 we take a convex combination of the previous and
the new participation matrices so as to obtain compatibility.  It is straightforward to prove that for any line $ij$, the quantity
$\bV((1 - t) \mA_{k-1} + t \hat \mA_k)_{ij}$ is a convex quadratic function over $0 \le t \le 1$ and so the stepsize computation can be performed exactly (see Appendix).

Typically, $\bmath{\Delta( \bar f^k, \hat s^2_k)} < \bmath{\Delta( \bar f^k, \hat s^2_{k-1})}$ (since $\bmath{(\hat s_k, \hat \mA_k)}$ solves the problem in Step 2). Thus, in order to
obtain a large improvement in variance metric, we want the stepsize $\lambda$ to be large.  This is
the reason why constraint \eqref{careful} is needed in problem
$\bmath{VShift(\bar f^k, \mA_{k-1}, \tau)}$: without such a constraint a line $ij$ in $\bT(\bar f^k, \mA^{k-1}, \tau)$ would enforce a short step if $\hat s_{k, ij}$ is large. Also, as stated the procedure may halt in Step 1. This need not be the case by relying on the expedient of resetting $\tau$ to a smaller value (e.g., half) and repeating Step 1 until feasibility is attained. A more comprehensive solution would be to perform a combination of Step 1 and Step 2, effectively a first-order gradient step, so as to move away from the current point.  We have not implemented this
patch as it did not prove necessary for small values of $\tau$ and $K$.

In the Appendix we will prove an important result concerning the above procedure under model (I) of the variance metric (Section \ref{varmetric}). 
Let $\Delta^*$ be the minimum variance metric over all compatible pairs, i.e.
\begin{align}
  \Delta^*  & \ \doteq \ \min \{ \, \Delta(\bar f, \bV(\mA)) \, : \, (\bar f, \mA) \, \mbox{compatible} \}, \label{mindelta}
\end{align}
\begin{theorem} \label{nonstop} Under model (I) if Procedure \eqref{varshift} stops at Step 5 of iteration $k$ then $\bmath{\Delta( \bar f^{k-1}, s^2_{k-1})} = \Delta^*.$
\end{theorem}
\section{Numerical examples}
Here we consider numerical examples based on  ``case2746wp'' available through MATPOWER \cite{matpowerpaper}, with $2746$ buses, $3514$ branches,  $520$ generators and sum of loads $24873$. We have added $22$ stochastic injection sites, with
sum of mean injections $4611.57$ (approx. $18.5 \%$ penetration) and average
ratio of standard deviation to mean of $0.3$.  Our experiment proceeded as follows.\\

\noindent {\bf Step 1.} First, solve the safety-constrained problem \eqref{scOPF1} with all safety parameters set to $3$ (three standard deviations). All generators are available for balancing.  This process required
29 iterations of the cutting-plane algorithm in \cite{ccopf} or \cite{line1} and
consumed approximately one CPU minute in a standard workstation, using CPLEX 12.6
\cite{cplex} as the underlying SOCP solver.\\

\noindent {\bf Step 2.} We then applied the Procedure \ref{varshift} with $\tau = 0.1$ and $K = 2$.  The set of buses used to balance stochastic deviations, $\mR$, was the subset of generating buses that was selected in Step 1, of cardinality $11$. We expand on Step 2 next.\\

We used, for variance metric, a nonconvex example of model (II) in Section \ref{varmetric} . The variance metric we used was
\begin{align}
  &  \sum_{ij \in \mF} s^2_{ij} \label{heredeltadef},
\end{align}
i.e. the sum of flow variances in lines in set $\mF$. At any iteration $k$ of procedure \ref{varshift}, the set $\mF$ used
to define \eqref{heredeltadef} is the union of two sets:
\begin{itemize}
\item [(a)] The $100$ lines with largest flow.
\item [(b)] Those lines for which the safety constraint is nearly binding, i.e.
  lines $ij$ for which $|\bar f_{ij}| + \nu_{ij} s_{ij} \ge (1 - \tau)f^{max}_{ij}$
  (``nearly-binding'' lines).
\end{itemize}
This problem is quintessentially non-convex and it is of interest to see whether
our iterative procedure can indeed reduce the metric.  A summary of the run is as follows:

\noindent {\bf Iteration k = 1, Step 1.} Requires $1.12$ seconds, optimal
expected generation cost $\approx 1.1\times10^{06}$.  \\
\noindent {\bf Iteration k = 1, Step 2.} The output of Step 1 produced a set of nearly-binding lines of cardinality $5$; so
that $|\mF| = 105$, with variance metric \eqref{heredeltadef} of value $6.3\times10^{04}$.  The optimization problem \eqref{scOPFshiftvar} had approximately $140000$ variables and a similar number of constraints and approximately one million
nonzeros. Its solution, using Gurobi 7.02 \cite{gurobi}, required $2.32$ seconds. The variances
$\bmath{\hat s^2_1}$ yield metric $\approx 2.3\times10^{04}$.\\
\noindent {\bf Iteration k = 1, Step 3.} Here, $\lambda \approx 0.55$.\\
\noindent {\bf Iteration k = 1, Steps 4 and 5.} The variances $\bmath{s^2_1}$ yield
a metric of value $\approx 4.65\times10^{04}$.\\ 
\noindent {\bf Summary.} One iteration of Procedure \ref{varshift} keeps expected
generation approximately constant and reduces variance metric by approximately $35 \%$.\\ \\
\noindent {\bf Iteration k = 2, Step 1.} Similar runtime, expected generation cost remains nearly same.\\
\noindent {\bf Iteration k = 2, Step 2.} The number of nearly-binding branches is now
$24$ with $\mF = 120$. Solution statistics are similar to those for iteration 1. The variance
metric attained by $\hat s^2_2$ is $2.89\times10^{04}$.\\
\noindent {\bf Iteration k = 2, Step 3.} Here, $\lambda \approx 0.29$.\\
\noindent {\bf Iteration k = 2, Steps 4 and 5.} The variances $\bmath{s^2_2}$ yield a metric of value $4.50\times10^{04}$.\\
\noindent {\bf Summary.} Two iterations of Procedure \ref{varshift} again
keep expected cost nearly constant, but reduce variance metric by roughly $40\%$
relative to the original value. \\

It is also instructive to look at the structure
of the power flows. At termination, the largest magnitude expected flow is of approximately of value $817$. The third largest flow has value $632$ and attains the largest single standard deviation, of value approximately $91$. In contrast, 
the $101^{st}$ largest flow value is approximately $144$ and lines at or below
that ranking of flow magnitude have much smaller standard deviation of flow; approximately $15$. 
Thus the procedure shifts variance away from high flow lines but also without
creating very high variance, low expected flow lines.  Of course, other choices
for the variance metric $\Delta$ will produce different tradeoffs.
\section{Conclusion}
In this work we have described efficient computational procedures that postprocess security-constrained DC-OPF solutions toward nearly-optimal solutions that attain
significantly lower variance metrics.  In future work we plan to explore alternative metrics, for example to capture engineering details, through
the use of derivative-free optimization \cite{nomad}.

\section{Acknowledgment}
This work was supported by the Advanced Grid Modeling Program
of the Office of Electricity at the U.S Department of Energy and
the Center for Nonlinear Studies at Los Alamos National Laboratory.

\bibliographystyle{IEEEtran}
\bibliography{apurv.bib}

\section{Appendix}
Here we provide the technical analysis underpinning Procedure \ref{varshift}.
\begin{remark} \label{quadrem} Let $\mA, \mA' \in \mK$.  Then for any
  $0 \le t \le 1$, $(1 - t)+ \mA +  t \mA' \in \mK$ and for any line $ij$, $\bV((1 - t)\mA +  t \mA' )_{ij}$ is a convex quadratic function of $t$.
\end{remark}
\noindent {\em Proof.} The first claim follows since $\mK$ is convex and the
second  using expression \eqref{Varexpression}. \qed
\begin{lemma} \label{shiftlemma1} Suppose that in iteration $k$ the algorithm reaches
  Step 2.  Let $(\bar f, \mA )$ be an arbitrary compatible pair.
Then there
exists $0 < \gamma \le 1$ such that for all $0 \le t \le \gamma$,
$$( \, (1 - t) \bar f^k + t \bar f \, , \, (1 -t) \mA_{k-1} + t \mA \,)$$ is a compatible pair. \end{lemma} 
\noindent \textit{Proof.}   Consider any line $ij$.  For real
$0 \le t \le 1$ let $s^2_{ij}(t) \ = \bV( ( 1 - t) \mA_{(k-1)} + t \mA)$.  By construction in Step 1,
$$ |\bar f^k_{ij}| + \nu_{ij} s_{ij}(0) \, \le \, (1 - \tau)f^{max}_{ij}.$$
Since $s^2_{ij}(t)$ is a quadratic function of $t$ and thus continuous, it follows
that there exists $\gamma_{ij} > 0$ such that
$$ |\bar f^k_{ij}| + \nu_{ij} s_{ij}(t) \, \le \, f^{max}_{ij}$$
for all $t \le \gamma_{ij}$. The quantity $\gamma_{ij}$ can be computed exactly using e.g. \eqref{Varexpression} to obtain an explicit representation of $s^2_{ij}(t)$. The result follows using $\gamma = \min_{ij} \gamma_{ij}$. \qed

Recall that $\Delta^*$ is the minimum variance metric over all compatible pairs
\eqref{mindelta}, and define
\begin{align}
(\bar f^*, \mA^*)  & \ \doteq \ \text{argmin} \{ \, \Delta(\bar f, \bV(\mA)) \, : \, (\bar f, \mA) \, \mbox{compatible} \}, \nonumber
\end{align}

\begin{lemma} \label{shiftlemma2} Suppose that in iteration $k$ the algorithm reaches   Step 2 and that  $\Delta^* < \Delta( \bar f^{k}, s_{k-1}^2)$. Let $\gamma$
  be as in Lemma \ref{shiftlemma1}, and write
  $$s^2(\gamma) \doteq \bV( ( 1 - \gamma) \mA_{(k-1)} + \gamma \mA^*).$$
Then
$$\Delta(\bar f^{k}, s^2(\gamma) ) \ < \ \Delta(\bar f^{k}, s^2_{k-1}).$$
\end{lemma}
\noindent \textit{Proof.} Apply Lemma \ref{shiftlemma1} with $\mA = \mA^*$
to obtain $\gamma$ as in (a). To prove (b) consider any line $ij$. Since 
for any line $ij$, the function $s^2_{ij}(t)$ is a convex
quadratic function of $t$,  $s^2_{ij}(\gamma) \le (1 - \gamma) s^2_{ij}(0) + \gamma s^2_{ij}(1)$, and
therefore, since $\Delta_{ij}(s_{ij}^2)$ is convex and nondecreasing in $s_{ij}$,
$$\Delta_{ij}(s^2_{ij}(\gamma)) \le (1-\gamma) \Delta_{ij}(s^2_{ij}(0)) +  \gamma \Delta_{ij}(s^2_{ij}(1)).$$
Summing this expression over $\mF(\bar f^k)$ we obtain
$$\Delta(\bar f^k, s^2(\gamma)) \le (1 - \gamma) \Delta(\bar f^k, s^2_{k-1}) + \gamma \Delta^* < \Delta(\bar f^k, s^2_{k-1}) $$
where the last inequality holds because $\gamma > 0$. \qed

{\bf Comment:} Lemma \eqref{shiftlemma2} shows that under very general conditions the
variance-shifting Step 2 of Procedure \ref{varshift} does lead to a reduction
in variance metric. In fact,

\begin{lemma} \label{shiftlemma3} Under model (I) of the variance metric
  $\Delta(\bar f^k, s^2_{k}) < \Delta(\bar f^{k-1}, s^2_{k-1})$ unless
  $\Delta(\bar f^{k-1}, s^2_{k-1}) = \Delta^*$.
 \end{lemma}
\noindent{\em Proof.} Under model (I) we have $\mF = \mE$ in \eqref{Deltadef}. Hence $\Delta(\bar f^k, s^2_{k-1}) = \Delta(\bar f^{k-1}, s^2_{k-1})$.  The result now follows from
    Lemma \ref{shiftlemma2}.
    In other words we obtain Theorem \ref{nonstop} .
\end{document}